\documentclass{amsart}

\usepackage{amscd}
\usepackage{amssymb}
\usepackage{url}
\input xy
\xyoption{all}

\newtheorem{anyprop}{Anyprop}[section]

\newtheorem{theorem}[anyprop]{Theorem}
\newtheorem{lemma}[anyprop]{Lemma}
\newtheorem{proposition}[anyprop]{Proposition}
\newtheorem{corollary}[anyprop]{Corollary}

\theoremstyle{definition}

\newtheorem{definition}[anyprop]{Definition}

\newtheorem{example}[anyprop]{Example}

\newtheorem{remark}[anyprop]{Remark}

\renewcommand{\AA}{\mathbb{A}}

\newcommand{\NN}{\mathbb{N}}
\newcommand{\ZZ}{\mathbb{Z}}

\newcommand{\FF}{\mathbb{F}}

\newcommand{\PP}{\mathbb{P}}

\newcommand  {\shA}     {\mathcal{A}}

\newcommand  {\shE}     {\mathcal{E}}
\newcommand  {\shF}     {\mathcal{F}}

\newcommand  {\shM}     {\mathcal{M}}

\newcommand  {\shN}     {\mathcal{N}}
\newcommand  {\shL}     {\mathcal{L}}

\newcommand  {\shX}     {\mathcal{X}}

\newcommand  {\fop}     {\mathfrak{p}}

\newcommand  {\Char}    {\operatorname{char}}

\newcommand  {\Ext}     {\operatorname{Ext}}

\newcommand  {\length}  {\operatorname{length}}

\newcommand  {\modu}     {\operatorname{mod}}

\renewcommand{\O}       {\mathcal{O}}

\newcommand  {\Pic}     {\operatorname{Pic}}

\newcommand  {\Proj}    {\operatorname{Proj}}

\newcommand  {\ra}      {\rightarrow}

\newcommand  {\rk}    {\operatorname{rk}}

\newcommand  {\Spec}    {\operatorname{Spec}}

\newcommand  {\Syz}     {\operatorname{Syz}}

\newcommand{\comdots}{ , \ldots , }
\newcommand{\komdots}{ , \ldots , }

\newcommand{\subsetdots}{ \subset \ldots \subset }

\newcommand{\lto}{\longrightarrow}

\newcommand{\expo}{{s}}
\newcommand{\expa}{{e}}
\newcommand{\scalar}{{\lambda}}

\newcommand{\degu}{{u}}

\newcommand{\algu}{{T}}
\newcommand{\alguval}{{t}}

\newcommand{\cloalg}[1]{{ \overline {#1 \!\!}}}

\theoremstyle{remark}

\numberwithin{equation}{section}

\begin{document}

\title[On deep Frobenius descent and flat bundles]
{On deep Frobenius descent and flat bundles}

\author[Holger Brenner and Almar Kaid]{Holger Brenner and Almar Kaid}

\address{Universit\"at Osnabr\"uck, Fachbereich 6: Mathematik/Informatik, Albrechtstr. 28a,
49069 Osnabr\"uck}

\email{hbrenner@uni-osnabrueck.de}

\address{Department of Pure Mathematics, University of Sheffield,
  Hicks Building, Hounsfield Road, Sheffield S3 7RH, United Kingdom}
\email{a.kaid@sheffield.ac.uk}

\thanks{}





\begin{abstract}
Let $R$ be an integral domain of finite type over $\ZZ$ and let
$f:\shX \ra \Spec R$ be a smooth projective morphism of relative
dimension $d \geq 1$. We investigate, for a vector bundle $\shE$
on the total space $\shX$, under what arithmetical properties of a
sequence $(\fop_n, e_n)_{n \in \NN}$, consisting of closed points
$\fop_n$ in $\Spec R$ and Frobenius descent data $\shE_{\fop_n}
\cong {F^{e_n}}^*(\shF)$ on the closed fibers $\shX_{\fop_n}$, the
bundle $\shE_0$ on the generic fiber $\shX_0$ is semistable.
\end{abstract}

\maketitle



\section{Introduction}

This paper is a continuation of our note \cite{brennerkaidfrobeniusdescent}.
Let $\ZZ \subseteq R$ be an integral domain of finite type and let
$f: \shX \to \Spec R$ be a smooth projective morphism of relative
dimension $d \geq 1$.
For a locally free sheaf $\shE$ on $\shX$ we are interested how
properties of $\shE_0=\shE|_{\shX_0}$ on
$\shX_0=\shX \times _{\Spec R}\Spec Q(R)$ (in characteristic zero)
are related to properties of $\shE_\fop$ on
$\shX_\fop= \shX_{\kappa(\fop)}= \shX \times _{\Spec R}\Spec \kappa(\fop)$
for closed points $\fop$ with finite residue class field (in positive
characteristic). More specifically, we study what Frobenius descent
properties for $\shE_\fop$ force $\shE_0$ to be semistable. We say
that $\shE_\fop$ admits an $e$th \emph{Frobenius descent} if
$\shE_\fop
\cong {F^e}^*(\shF)$ for some locally free sheaf $\shF$ on
$\shX_\fop$, $F$ being the absolute Frobenius morphism.
In \cite{brennerkaidfrobeniusdescent} we gave, addressing a question of K. Joshi \cite{joshiremarkonvectorbundles},
an example of a bundle
on a smooth projective relative curve over $\Spec \ZZ$ such that
$\shE_p$ admits a first Frobenius descent ($e=1$) for infinitely
many prime numbers $p$, yet $\shE_0$ is not semistable.

Here we allow the Frobenius exponents $e$ to grow. We will consider
\emph{Frobenius descent sequences}
$(\fop_n, e_n)_{n \in \NN}$ consisting of closed points $\fop_n$ and natural numbers $e_n$
such that $\shE_{\fop_n}$ admits an $e_n$th Frobenius descent.
It turns out for $\shX$ being a relative curve that
one has to compare the growth of $e_n$ with the number of elements
in the residue class fields $\kappa(\fop_n)$. Our main positive result, Theorem
\ref{maintheorem}, states that $\shE_0$ is semistable if
$e_n - |\kappa(\fop_n)|^c \to \infty$, where $c$ is a constant depending on the
genus of the curve and the rank of the bundle. We will give examples
of Frobenius descent sequences where the characteristic is constant
and
$e_n
\to \infty$, but $\shE_0$ is not semistable (Examples \ref{counterexampleconstantchar} and \ref{counterexampleconstantthree}).

We give a quick overview for the organization of this paper:
In Section \ref{remarksflatfrobenius} we deal with bundles over
varieties over a field of fixed positive characteristic. We recall the notions of
strong semistability,  flat bundles (introduced by Gieseker),
Frobenius descent and Frobenius periodicity.
We exhibit how these concepts are related over a finite field and how they differ in general.
In Section \ref{deeperfrobrelativecurves} we prove after some preparatory work our main theorem,
and Section \ref{counterexamplesection}
presents the examples.

\section{Flat bundles, Frobenius descent and strong semistability}
\label{remarksflatfrobenius}

To begin with we recall some essential notions regarding
semistable sheaves on a smooth projective variety $X$ over an algebraically closed field $K$ (for a
complete account see the book \cite{huybrechtslehn}). A coherent
torsion-free sheaf $\shE$ on $X$ is \emph{semistable} (in the
sense of Mumford and Takemoto) if for every coherent subsheaf $0
\neq \shF \subset \shE$ the inequality $\mu(\shF) =
\deg(\shF)/\rk(\shF) \leq \deg(\shE)/\rk(\shE) = \mu(\shE)$ holds.
The sheaf $\shE$ is \emph{stable} if the inequality is always
strict. The degree of a sheaf $\shF$ is defined using intersection
theory and a fixed ample invertible sheaf $\O_X(1)$ as $\deg
(\shF) = \deg (c_1(\shF) . \O_X(1)^{\dim (X)-1})$. For every
coherent torsion-free sheaf $\shE$ there exists the
\emph{Harder-Narasimhan filtration} $\shE_1 \subset \shE_2
\subsetdots \shE_t =\shE$ such that $\shE_i/\shE_{i-1}$ is
semistable and $\mu(\shE_1) > \mu(\shE_2/\shE_1) > \ldots >  \mu
(\shE/\shE_{t-1})$. The slopes $\mu(\shE_1)$ and $\mu
(\shE/\shE_{t-1})$ are also denoted
by $\mu_{\max}(\shE)$ and $\mu_{\min}(\shE)$ respectively. If $K$ is not algebraically closed, then we define the 
terms degree, semistable, etc. via the algebraic closure of $K$.

On a smooth projective variety $X$ over a field $K$ of positive
characteristic $p > 0$ we consider the absolute Frobenius morphism
$F: X \ra X$ which is the identity on the topological space $X$
and the $p$th power map $(-)^p$ on the structure sheaf $\O_X$.
A coherent torsion-free sheaf $\shE$ is called \emph{strongly
semistable} if all Frobenius pull-backs ${F^e}^*(\shE)$ are
semistable. The Frobenius pull-back of an algebraic cycle $Z_i \in
A^{i}(X)$ of codimension $i$ is $p^i Z_i$. Hence the degree of
$F^*(\shE)$ with respect to $\O_X(1)$ is
$\deg( c_1(F^*(\shE)) . \O_X(1)^{\dim (X) -1}) = p \deg (\shF) $.

\begin{definition}
Let $X$ be a smooth projective variety over a field of positive characteristic and let $\shE$ be a locally
free sheaf on $X$. We say that $\shE$ admits an \emph{$e$th
Frobenius descent} if there exists a locally free sheaf $\shF$
such that
$\shE \cong {F^e}^*(\shF)$. We say that $\shE$ admits an \emph{infinite
Frobenius descent} if $\shE \cong {F^e}^*(\shF_e)$ holds for every
$e \in \NN$ and locally free sheaves $\shF_e$. Furthermore, the
locally free sheaf $\shE$ is called \emph{flat}, if there exists a
collection $(\shF_e)_{e \in \NN}$ of locally free sheaves such that
$\shF_0 = \shE$ and $\shF_e \cong F^*(\shF_{e+1})$.
\end{definition}

The notion of flat (sometimes called stratified) bundles was
introduced by Gieseker (cf. \cite[Definition
1.1]{giesekerflatfundamental}). We emphasize  that the notion of flatness depends on the 
base field (see Example \ref{linebundleexample}). If in positive characteristic a
semistable vector bundle $\shE$ admits a Frobenius descent $\shE
\cong F^*(\shF)$ then clearly $\shF$ has to be semistable too. The
following easy observation is well-known (cf. \cite[Lemma 5]{giesekerfrobenius}), but for completeness  we present it with a proof.

\begin{lemma}
\label{mumaxboundsemistable}
Let $X$ be a smooth projective variety defined over a field of
positive characteristic $p > 0$, let $\O_X(1)$ be a fixed ample line
bundle and let
$\shE$ be a locally free sheaf of degree zero and rank $r$ on $X$
with
$b:=\mu_{\max}(\shE)$. If $\shE \cong {F^e}^*(\shF)$ holds for some
$e > \log_2(rb)$ and a locally free sheaf $\shF$ then $\shF$ is
semistable. Moreover, if $\shE $ is a flat vector bundle on
$X$ then there exists an $e_0$ such that all descent bundles
$\shF_e$ are semistable {\rm(}with respect to every $\O_X(1)${\rm)} for $e \geq e_0$.
\end{lemma}

\begin{proof}
Assume $\shF$ is not semistable. Since by assumption $\shE$ has
degree zero and rank $r$, so has $\shF$. Therefore, we have
$\mu_{\max}(\shF) \geq 1/r$ and we obtain (we may assume that $\dim
(X) \geq 1$)
$$\mu_{\max}({F^e}^*(\shF)) \geq \frac{p^e}{r} \geq \frac{2^e}{r}> \frac{2^{\log_2 (rb)}}{r} = b.$$
This yields a contradiction. For a flat vector bundle $\shE$ we have
by definition
$\shF_e \cong F^*(\shF_{e+1})$ for every $e \in \NN$. Hence,
the supplement follows now from the equalities $\shE = \shF_0 =
F^*(\shF_1) = F^*(F^*(\shF_2)) = {F^2}^*(\shF_2) = \ldots =
{F^e}^*(\shF_e)$ for every $e>0$.
\end{proof}

\newcommand{\etcov}{{\psi}}

\begin{proposition}\label{equivalentcondiontsflat}
Let $X$ be a smooth projective variety defined over a finite field
$K$ and let $\shE$ be a locally free sheaf of rank $r$ on $X$.
Then the following conditions are equivalent:
\begin{enumerate}
\item $\shE \cong {F^e}^*(\shE)$ for some $e >0$.
\item There exists
an $\acute{e}$tale cover $\etcov:\tilde{X} \ra X$ such that
$\etcov^*(\shE)
\cong \O_{\tilde{X}}^r$.
\item $\shE$ is a flat vector bundle.
\item $\shE$ admits infinite Frobenius descent.
\end{enumerate}
\end{proposition}

\begin{proof}
\noindent $(1) \Leftrightarrow (2)$. This follows from \cite[Satz
1.4]{langestuhler}; see also Remark \ref{langestuhlerremark}.

\noindent $(1) \Rightarrow (3)$. We define
$\shF_i:={F^{(-i \modu e)}}^*(\shE)$ (taking $0\leq -i \modu e < e $) and get
\[F^*(\shF_{i+1}) = F^* ({F^{((-i-1) \modu e)}}^*(\shE) )
= \begin{cases}  {F^e}^*(\shE) \cong \shE \cong \shF_i \mbox{ if } i \equiv 0 \!\!\!\! \mod e,\\
{F^{(-i \modu e)}}^*(\shE) \cong \shF_i \mbox{ if } i \not
\equiv 0 \modu e
\end{cases} \]
and see that $\shE$ has the structure of a flat vector bundle.

\noindent $(3) \Rightarrow (4)$. As mentioned in the proof of
Lemma \ref{mumaxboundsemistable} we have $\shE  ={F^e}^*(\shF_e)$ for every $e > 0$ which yields infinite
Frobenius descent.

\noindent $(4) \Rightarrow (1)$. By assumption there exist locally free sheaves $\shF_e$, $e \in \NN$, such that $\shE \cong
{F^e}^* (\shF_e)$ for every $e > 0$. This implies  the equalities
$c_i(\shE) = p^{ei} c_i(\shF_e)$ for the
Chern classes of $\shE$ for all $e$, $p = \Char (K)$, $i=1 \komdots \dim X$. Hence, all the
Chern classes $c_i(\shE)$ and $c_i(\shF_e)$ are numerically
trivial for $i=1 \komdots \dim X$ and $e \geq 0$, because
$ \deg(c_i(\shE). Z) =
\deg (p^{ei} c_i(\shF_e). Z)=p^{ei} \deg(c_i(\shF).Z)$ for every algebraic cycle $Z$ of dimension $i$.
By
Lemma
\ref{mumaxboundsemistable} the bundles $\shF_e$ are semistable for
$e$ sufficiently large. Since the family of semistable vector
bundles of fixed rank with numerically trivial Chern classes is bounded by
Langer's theorem \cite[Theorem 4.2]{langersemistable} and since we are working over
a finite field, there are only finitely many isomorphism classes of
such bundles. Hence we must have a repetition $\shF_t \cong \shF_s$
for $t>s$. Applying ${F^t}^*$ yields
$$\shE \cong {F^t}^*(\shF_t) \cong {F^t}^*(\shF_s) \cong {F^{(t-s)}}^*({F^s}^*(\shF_s)) \cong
{F^{(t-s)}}^*(\shE)$$ and we obtain a repetition of the vector
bundle
$\shE$.
\end{proof}

In \cite[Theorem 1.10]{giesekerflatfundamental} D. Gieseker proved
that if every flat bundle on a smooth projective variety $X$ over an
algebraically closed field $K$ is trivial then the algebraic
fundamental group
$\pi_1(X)$ is trivial. It is an open question in
\cite[Beginning of \S 2]{giesekerflatfundamental} whether the
converse is true (and was proven for projective spaces and
$K3$ surfaces in \cite[Theorem 2.2 and
Theorem 2.3]{giesekerflatfundamental}). For bundles which are
flat over a finite field we have the following converse.

\begin{corollary}
Let $X$ be a smooth projective variety defined over a finite field
$K$ with $\pi_1(X_{\bar{K}})=0$. Then every flat bundle on $X$ is trivial.
\end{corollary}
\begin{proof}
Let $\shE$ be a flat bundle on $X$. By Proposition
\ref{equivalentcondiontsflat} there exists an $\acute{e}$tale
cover $\etcov: \tilde{X} \ra X$ such that $\etcov^*(\shE)$ becomes
trivial on $\tilde{X}$. Hence $\etcov^*(\shE)_{\bar{K}}$ is trivial
on
$\tilde{X}_{\bar{K}}$. By assumption
$\pi_1(X_{\bar{K}})$ is trivial and therefore $X_{\bar{K}}$ is simply
connected (cf. \cite[Theorem 5.3]{milne}). Hence $\shE_{\bar{K}}$ is
trivial on
$X_{\bar{K}}$ and so $\shE$ is trivial on $X$.
\end{proof}

We also present a weaker version of Proposition
\ref{equivalentcondiontsflat} which goes essentially back to
\cite{langestuhler}.

\begin{proposition}
\label{equivalentcondiontsstronglysemistable}
Let $X$ be a smooth projective variety defined over a finite field
$K$ and let $\shE$ be a locally free sheaf. Then the following
conditions are equivalent.
\begin{enumerate}
\item There exists a repetition ${F^t}^*(\shE) \cong F^{s*}
(\shE)$ with $t > s$. \item $\shE$ is strongly semistable  $($with
respect to every $\O_X(1))$ with numerically trivial Chern
classes. \item There exists an $\acute{e}$tale cover
$\etcov:\tilde{X} \ra X$ such that $\etcov^*({F^s}^*(\shE))$
becomes trivial for a certain $s>0$.
\end{enumerate}
\end{proposition}

\begin{proof}
\noindent $(1)\Rightarrow(2)$. The repetition yields the equations
$p^{t i} c_i(\shE) = p^{s i} c_i(\shE)$ for $i=1 \comdots \dim (X)$,
hence $p^{(t-s) i} c_i(\shE) = c_i(\shE)$ and the assertion on the
Chern classes follows. In particular, $\deg(\shE)=0$ for every $\O_X(1)$. We have $F^{t*}(\shE) =
F^{(t-s)*}({F^s}^*(\shE)) \cong {F^s}^*(\shE)$ and so $F^{k(t-s)*}({F^s}^*(\shE)) \cong {F^s}^*(\shE) $ for all $k$.
Thus ${F^s}^*(\shE)$ is semistable by Lemma
\ref{mumaxboundsemistable}. Every ${F^e}^*(\shE)$ has a Frobenius
pull-back which is isomorphic to $ {F^s}^*(\shE) $, hence $\shE$ is
strongly semistable.

\noindent $(2)\Rightarrow(1)$. This follows also immediately from
the boundedness of the family of semistable vector bundles with
numerically trivial Chern classes.

\noindent $(1)\Leftrightarrow(3)$. For this equivalence see
\cite[Satz 1.4]{langestuhler}.
\end{proof}

\begin{corollary}
Let $X$ be a smooth projective variety defined over a finite field
$K$ and let $\shE$ be a flat vector bundle on $X$.
Then $\shE$ is strongly semistable.
\end{corollary}

\begin{proof}
This follows from Proposition \ref{equivalentcondiontsflat} and Proposition \ref{equivalentcondiontsstronglysemistable}.
\end{proof}

\begin{remark}
\label{langestuhlerremark}
The implications $(1)\Rightarrow (2) \Rightarrow (3) \Rightarrow
(4)$ of Proposition \ref{equivalentcondiontsflat} and $(1)
\Rightarrow (3) \Rightarrow (2)$ of Proposition
\ref{equivalentcondiontsstronglysemistable} also hold for smooth
projective varieties defined over an arbitrary field of positive
characteristic. The cited result \cite[Satz 1.4]{langestuhler} of Lange and Stuhler holds
over a finite field, but not in the generality stated there. An example of Y. Laszlo mentioned in \cite[after Theorem 1.1]{biswasducrohetlangestuhler} shows that for a semistable $\acute{e}$tale trivializable vector bundle $\shE$ we do not necessarily have a Frobenius periodicity ${F^e}^*(\shE) \cong \shE$ for an $e>0$. We provide also such an example in \ref{antilangestuhlerexample}. It was also shown in \cite[Theorem 1.1]{biswasducrohetlangestuhler} that the theorem of Lange and Stuhler is true if $\shE$ is stable.

A flat vector bundle $\shE $ over an infinite field
of positive characteristic is in general not semistable. Examples
for this behavior were constructed in
\cite{giesekerfrobenius} by D. Gieseker. 
\end{remark}

In the following we give examples
of flat and strongly semistable bundles without any
Frobenius repetitions.

\begin{example}
\label{linebundleexample}
Let $K$ be an algebraically closed field of positive characteristic
$p$ and let $X$ be a smooth projective variety defined over $K$. If
$\shL \in \Pic^0(X)$ is a line bundle algebraically equivalent to zero then $\shL$ is
flat since $\Pic^0(X)$ is $p$-divisible; see
\cite[Application II.6.2]{mumfordabelian}.
Suppose furthermore that
$K$ contains a transcendental element and $\shL$ is non-torsion in
$\Pic^0(X)$. Then all powers ${F^e}^*(\shL) \cong \shL^{p^e}$ are distinct for $e \geq
0$. Hence, we obtain a flat vector bundle of rank one without any
repetitions. By \cite[Theorem 1.1]{biswasducrohetlangestuhler} $\shL$ is not 
$\acute{e}$tale trivializable.

In contrast, if $X$ is defined over a finite field  $K$ such that
the $p$-torsion of $\Pic^0(X)$ is non-trivial, we find a line bundle
$\shL \neq \O_X$ such that $\shL^p \cong \O_X$. Hence, there is no repetition
$\shL \cong \shL^{p^{e}}$ for any $e >0$.
By Proposition \ref{equivalentcondiontsflat} the line bundle
$\shL$ can not have infinite Frobenius descent. In particular, the property of being flat depends on the base field.
\end{example}

The following Proposition gives a method to construct flat vector
bundles of rank two without any Frobenius repetition. Unlike the previous example
these bundles are $\acute{e}$tale trivializable and hence come by \cite[Proposition 1.2]{langestuhler} from a representation $\rho:
\pi_1(X) \ra GL_2(K)$ of the algebraic fundamental group $\pi_1(X)$.

\begin{proposition}
\label{flatbundleconstruction}
Let $X$ be a smooth projective variety defined over an
algebraically closed field $K$ of characteristic $p > 0$ which
contains transcendental elements. If the $p$-rank of $X$ is at least
two then there exist non-trivial extensions of $\O_X$ by $\O_X$
which are flat, strongly semistable and without any Frobenius repetitions.
\end{proposition}

\begin{proof}
The Frobenius morphism acts on $V:=H^1(X,\O_X)$ $p$-linearly, i.e.,
$F^*(\lambda c) = \lambda^p F^*(c)$ for all $\lambda \in K$
and $c \in V$. We decompose $V=V_s \oplus V_n$ in its semi-simple
part $V_s$, where there exists a basis $v_1 \komdots v_k$ such that $F^*(v_i)=v_i$ for all $i=1\komdots k$, and its nilpotent part
$V_n$, where $F^*$ is a nilpotent map (cf. \cite[III.14, last
corollary]{mumfordabelian}). The assumption on the $p$-rank of
$X$ means that $\dim_K V_s \geq 2$. We fix two linearly independent
vectors $v$ and $w$ in $V_s$ with $F^*(v)=v$ and $F^*(w)=w$, and we
consider the cohomology class
$c=v+tw$, where $t \in K^\times$ is a transcendental element.

Let the vector bundle $\shE$ be the extension of $\O_X$ by $\O_X$ given by $c \in
H^1(X,\O_X)= \Ext^1(\O_X,\O_X)$. Such a bundle is strongly
semistable. The $e$th Frobenius pull-back of $\shE$ is the extension
given by ${F^e}^*(c) = v+ t^{p^e}w$. The same reasoning shows that
$\shE$ itself is the $n$th Frobenius pull-back of the extension given by $v+ t^{1/p^n}
w$ (where $t^{1/p^n}$ exists since $K$ is algebraically closed).
Hence $\shE$ is flat.
Assume ${F^\expo}^*(\shE) \cong {F^\expa}^*(\shE)$ with $\expo >
\expa$. Then the corresponding cohomology classes must differ by a
scalar, say $ v+ t^{p^\expo}w= \scalar ( v+ t^{p^\expa}w)$.
Hence $\scalar=1$ and so $t^{p^\expo}=t^{p^\expa}$ and $t^{p^{\expo
-\expa}}=1$, contradicting the transcendence of $t$.
\end{proof}

\begin{example}
\label{antilangestuhlerexample}
We consider the Fermat quartic
$$C:=\Proj(\overline{\FF_5(t)}[X,Y,Z]/(X^4+Y^4-Z^4))$$ over the
algebraic closure of the function field $\FF_5(t)$. The
$\rm\check{C}$ech-cohomology classes $\frac{\sqrt[4]{2}
Z^3}{X^2Y}$ and $\frac{\sqrt[4]{2} Z^3}{XY^2}$ are linearly
independent in $V:=H^1(C,\O_C)$. The Frobenius acts on
$\frac{\sqrt[4]{2} Z^3}{X^2Y}$ as follows:
\begin{eqnarray*}
F^* \left ( \frac{\sqrt[4]{2} Z^3}{X^2Y} \right) \!\!& =&  \!\!\frac{2
\sqrt[4]{2} Z^{15}}{X^{10}Y^5} \,\,=\,\, \frac{2 \sqrt[4]{2} Z^3
Z^{12}}{X^{10}Y^5}\\  \!\!&=& \!\!
\frac{2 \sqrt[4]{2} Z^3 (X^{12}+3X^8Y^4+3X^3Y^8+Y^{12})}{X^{10}Y^5}\\
 \!\!&=&  \!\!\frac{\sqrt[4]{2} Z^3 (X^8Y^4)}{X^{10}Y^5} \,\,=\,\,
\frac{\sqrt[4]{2} Z^3}{X^2Y}
\end{eqnarray*}
and analogously we obtain $F^*(\frac{\sqrt[4]{2}
Z^3}{XY^2})=\frac{\sqrt[4]{2} Z^3}{XY^2}$  by symmetry, so
$\frac{\sqrt[4]{2} Z^3}{X^2Y}$ and $\frac{\sqrt[4]{2} Z^3}{XY^2}$
are classes which are fixed under the Frobenius. Hence the bundle defined by  $c:=\frac{\sqrt[4]{2}
Z^3}{X^2Y}+ t\frac{\sqrt[4]{2} Z^3}{XY^2} $
has all properties described in Proposition
\ref{flatbundleconstruction}.
The descent data are given by the extensions defined by
$c_n:=\frac{\sqrt[4]{2}
Z^3}{X^2Y}+\sqrt[5^n]{t}\frac{\sqrt[4]{2} Z^3}{XY^2}$, $n \geq 0$.
\end{example}

The following example was used in
\cite{brennermonskytightclosure} to show that tight closure does not commute with localization.
It also shows that strong semistability is not an open property for
geometric deformations (see also \cite[Example 2.11]{langersurvey} for another instructive example), though the existence of a Frobenius
repetition is. We will come back to this example in
Lemma
\ref{pullbackhnfiltration} and
Example \ref{counterexampleconstantchar}.
For a homogeneous coordinate ring $R$
of a projective variety $X$ and homogeneous elements
$f_1 \komdots f_n \in R$ of degrees
$d_i:= \deg(f_i)$, $i=1 \komdots n$,
such that the ideal generated by the $f_i$ is $R_+$-primary, the
\emph{syzygy bundle} $\Syz(f_1 \komdots f_n)(m)$ is given by the
short exact sequence
$$0 \lto \Syz(f_1 \komdots f_n)(m)\lto \bigoplus_{i=1}^n
\O_X(m-d_i) \stackrel{f_1 \komdots f_n}\lto \O_X(m) \lto 0.$$
It has rank $n-1$ and degree $ ( (n-1)m -\sum_{i=1}^n d_i ) \deg
(\O_X(1))$.

\begin{example}
We consider the smooth plane curve
$$C:=\Proj \,(\FF_2(\algu)[X,Y,Z]/(Z^4+XYZ^2+X^3Z+Y^3Z+(\algu +\algu^2) X^2Y^2))\, ,$$
where $\algu$ might be algebraic or transcendental. We consider the
syzygy bundle $\Syz(X^2,Y^2,Z^2)(3)$ and its
Frobenius pull-backs
$$\Syz(X^{2q},Y^{2q},Z^{2q})(3q) \cong
{F^e}^*(\Syz(X^2,Y^2,Z^2)(3)),$$ $q=2^e$. For $\algu$ transcendental
all these pull-backs are semistable, but this is false for every
algebraic value; see \cite[Theorem 4.13]{monskypoints4quartics} and
\cite[Corollary 4.6]{brennerhilbertkunz}. In the transcendental case
there are however no repetitions in the sequence
${F^e}^*(\Syz(X^2,Y^2,Z^2)(3))$. If there would be such a
repetition then this would also hold for almost all algebraic
values, contradicting the fact that in the algebraic case the
syzygy bundle $\Syz(X^2,Y^2,Z^2)(3)$ is not strongly semistable.
We do not know whether the bundle $\Syz(X^2,Y^2,Z^2)(3)$ is flat for
$\algu$ transcendental.
\end{example}

\begin{example}
Let $C$ be a smooth projective curve of genus two defined over an algebraically
closed field $K$ of characteristic $3$. It was shown in
\cite[Corollary 6.6]{laszlopaulyfrobeniuskummer} that the rational
map
$$V: \shM_C(2,\O_C) \dashrightarrow \shM_C(2,\O_C),~[\shE]
\longmapsto [F^*(\shE)]$$ (also called the \emph{Verschiebung})
from the moduli space $\shM_C(2,\O_C)$ parameterizing semistable
rank two vector bundles with trivial determinant is surjective.
This means for a stable rank two vector bundle $\shE$ that $\shE
\cong F^*(\shF)$, where $\shF$ is necessarily stable too. Hence
every stable rank two vector bundle is flat.

For a semistable but not stable bundle $\shE$ we obtain just sequences $[\shE] =
[\shF_0]=[F^*(\shF_1)]=[{F^2}^*(\shF_2)]= \ldots =
[{F^n}^*(\shF_n)]$ of $S$-\emph{equivalence classes} for every $n
\geq 0$ such that $[\shF_n] = [F^*(\shF_{n+1})]$. Such an
\emph{$S$-flat bundle} is in general not flat. For instance, let $C$
be an elliptic curve and consider the bundle $F_2$ in Atiyah's
classification
\cite{atiyahelliptic}, i.e., the unique indecomposable sheaf of rank
two and degree zero with
$\Gamma(C,F_2) \cong K$. The bundle $F_2$ is sitting inside the
short exact sequence
$$0 \lto \O_C \lto F_2 \lto \O_C \lto 0$$
and we see that $[F_2]=[\O_C^2]$. In particular, $F_2$ is
$S$-flat. If the Hasse invariant of $C$ is zero, we have $F^*(F_2)
\cong \O_C^2$. Therefore, $F_2$ is not flat since it is not a
pull-back of another vector bundle. If the Hasse invariant of $C$
equals one, then $F^*(F_2) \cong F_2$ and $F_2$ is flat.
\end{example}

\section{Deep Frobenius descent on relative curves}
\label{deeperfrobrelativecurves}

In this section we consider a smooth projective morphism $f: \shX
\ra \Spec R$ of schemes of relative dimension $d \geq 1$ where the ring $R$
is Noetherian. In this situation
we denote by $\shX_\fop:= \shX \times_{\Spec R} \Spec
\kappa(\fop)$ the fiber over the prime ideal $\fop \in \Spec R$,
where $\kappa(\fop):=R_{\fop}/\fop R_\fop$ denotes the residue class field
at the point $\fop$. These fibers are smooth projective schemes of
dimension $d$. Further, we fix an $f$-very ample line bundle
$\O_\shX(1)$ on $\shX$, in particular $\O_{\shX}(1)|_{\shX_{\fop}}$ is a very
ample line bundle on $\shX_{\fop}$. If $\shE$ is a locally free
sheaf on $\shX$ then $\shE_{\fop}:=\shE|_{\shX_\fop}$ denotes the
restriction of $\shE$ to the fiber $\shX_\fop$.

The following Lemma is well-known, but since it is essential for
our further progress we will present it with a full proof.

\begin{lemma}
\label{mumaxbound}
Let $R$ be a Noetherian ring and let $f: \shX \rightarrow \Spec R$ be 
a smooth projective morphism of relative dimension $d \geq 1$ together 
with a fixed $f$-very ample line bundle $\O_\shX(1)$. If $\shE$ is a locally free sheaf on
$\shX$ then there exists a global bound $b$ such that
$\mu_{\max}(\shE_{\fop}) \leq b$ for all $\fop \in \Spec R$.
\end{lemma}

\begin{proof}
Without loss of generality, we can assume that $\Spec R$ is connected. The locally free sheaf $\shE(m)$ is globally generated for $m \gg 0$ (cf. \cite[Theorem 5.17]{hartshornealgebraic}), i.e.,
we have a surjective morphism
$$\O_{\shX}^s(-m) \lto \shE \lto 0$$
for some $s > 0$. Restriction to each fiber $\shX_\fop$, $\fop \in \Spec R$, yields
$$\O_{\shX_\fop}^s(-m) \lto \shE_\fop \lto 0.$$
Hence $\mu_{\min}(\shE_\fop) \geq -m \deg(\O_\shX(1))$. In particular, this bound holds for the slope of every quotient sheaf of $\shE_\fop$. Therefore, if $\shE_\fop$ is not semistable and $\shF$ is the maximal destabilizing subsheaf of $\shE_\fop$ we obtain the simple inequalities

\begin{eqnarray*}
\mu_{\max}(\shE_\fop)&=&\mu(\shF) = \frac{\deg(\shF)}{\rk(\shF)}=\frac{\deg(\shE_\fop)}{\rk(\shF)} - \frac{\deg(\shE_\fop/\shF)}{\rk(\shF)}\\ &=& \frac{\deg(\shE)}{\rk(\shF)} - \frac{\rk(\shE_\fop/\shF)}{\rk(\shF)} \mu(\shE_\fop/\shF) \leq  \frac{\deg(\shE)}{\rk(\shF)} +
\frac{\rk(\shE_\fop/\shF)}{\rk(\shF)} m \deg(\O_\shX(1))\\ &\leq& \max(0, \deg(\shE)) + \rk(\shE) m \deg(\O_\shX(1)),
\end{eqnarray*}
where $\deg(\shE)=\deg(\shE_\fop)$ is constant in the family (cf. \cite[Definition 1.2.11]{huybrechtslehn} and \cite[Section II.5]{mumfordabelian}). This bound also holds for $\shE_\fop$ semistable.
\end{proof}

\begin{lemma}
\label{bundlebound}
Let $C$ be a smooth projective curve defined over a finite field
$K$. Then the number of isomorphism classes of semistable vector
bundles of rank $r$ and degree zero is bounded from above by
$|K|^c$, where $c$ is a constant depending only on $r$ and the genus of
$C$.
\end{lemma}

\begin{proof}
Let $\shE$ be a semistable vector bundle on $C$ of rank $r$ and
degree zero. Further, we fix an ample line bundle $\O(1)$ of degree
one. For $\ell>2g-1$  we have $\mu(\shE(\ell))=\ell > 2g-1$. Therefore, by \cite[Lemme 20]{seshadrifibre}, the bundle $\shE(\ell)$ is generated by global sections and $h^1(C,\shE(\ell))=0$. So we obtain a short exact sequence
$$0 \lto \shF\lto \O(-\ell)^s \lto \shE \lto 0$$
where $s:=h^0(C,\shE(\ell))=r \ell+r(1-g)$ and $\shF$ is a locally free sheaf of rank
$s-r$ and degree $-\ell s$. Let $0=\shF_0 \subset \shF_1 \subset \ldots \subset \shF_t = \shF$
be the Harder-Narasimhan filtration of $\shF$. Since $\mu(\shF_i) \leq -\ell$ holds for all $i=1 \komdots t$, the degrees of the subbundles $\shF_i$ are all negative. Thus we have $\deg(\shF/\shF_{t-1})\geq \deg(\shF)= -\ell s$ and moreover $\mu(\shF/\shF_{t-1})\geq -\ell s$. This yields
\begin{eqnarray*}
\mu((\shF/\shF_{t-1})(\ell(s+1))) &=& \mu(\shF/\shF_{t-1}) + \ell(s+1)\\
&\geq&-\ell s + \ell(s+1)\\
&=& \ell.
\end{eqnarray*}
So for $k:=\ell(s+1)$, again by \cite[Lemme 20]{seshadrifibre}, all twisted quotients $(\shF_i/\shF_{i-1})(k)$, $i=1 \komdots t$, are globally
generated with $h^1(C,(\shF_i/\shF_{i-1})(k))=0$.

We show now by induction that the bundles $\shF_i(k)$, and in particular $\shF(k)$, are also globally generated with $h^1(C,\shF_i(k))=0$. To do that, we look at the short exact sequences
$$0 \lto \shF_{i-1}(k) \lto \shF_i(k) \lto (\shF_{i}/\shF_{i-1})(k) \lto 0.$$
By the induction hypothesis $\shF_{i-1}(k)$ is generated by its global sections and we have
$h^1(C,\shF_{i-1}(k))=0$ which implies $h^1(C,\shF_{i}(k))=0$ too, since the quotients do not have
higher cohomology. Hence, we obtain a commutative diagram
$$ \xymatrix @-5pt {0 \ar[r] & \O \otimes \Gamma(C,\shF_{i-1}(k)) \ar[d] \ar[r] & \O \otimes \Gamma(C,\shF_i(k)) \ar[d] \ar[r] & \O\otimes \Gamma(C,(\shF_i/\shF_{i-1})(k)) \ar[d] \ar[r] & 0\\0  \ar[r]  & \shF_{i-1}(k) \ar[d] \ar[r] & \shF_{i}(k) \ar[r] & (\shF_i/\shF_{i-1})(k) \ar[d] \ar[r] & 0\\& 0          &            & 0        &    }
              $$
with exact rows and columns. We have to show that the morphism of sheaves in the middle is surjective too. But, since surjectivity is a local property, we can consider this diagram locally at a point $P \in C$. Then the surjectivity follows from the five lemma (cf. \cite[\S 42, Aufgabe 14]{schejastorch1}). So we have proved that the bundle $\shF(k)$ is globally generated and we have the following commutative diagram
$$ \xymatrix @-5pt {& \O(-k)^m \ar[d]\ar[dr]^\shA            &             &        & \\
0\ar[r]  & \shF \ar[d] \ar[r] & \O(-\ell)^{s} \ar[r]  & \shE \ar[r]  & 0\\
             & 0,                  &                      &                   &}
              $$
where $m:=h^0(C,\shF(k))=-\ell s + (s-r)k+(s-r)(1-g)$ and $\shA$ is an $m \times s$ matrix
with entries in $\Gamma(C,\O(k-\ell))=\Gamma(C,\O(\ell s))$. So every semistable
vector bundle of rank $r$ and degree zero is parameterized by such a
matrix $\shA$. The $K$-vector space
$\Gamma(C,\O(\ell s))$ has by Riemann-Roch dimension
$h^0(C,\O(\ell))=\ell s+1-g=:n$. Hence $|K|^c$ , with $c:=nms$,
is an upper bound for the number of possible matrices and we see
that the constant $c$ only depends on $r$ and
$g$.
\end{proof}

\begin{remark}
We used essentially Seshadri's result \cite[Lemme 20]{seshadrifibre}
of the existence of a common bound $m$ such that all semistable
vector bundles on the curve $C$ of rank $r$ and degree $\ell> m$ are
generated by their global sections. It is shown in \cite[Example 3.2]{langersurvey} that one can not expect such a bound for higher dimensional varieties if one fixes only the first Chern class of the bundles. Nevertheless, it is in principal
possible to generalize Lemma \ref{bundlebound} and the following theorem for semistable vector bundles of fixed rank and numerically trivial Chern classes to arbitrary dimension using \cite[Theorem 3.4]{langersurvey}.
\end{remark}

We are now able to prove our main theorem. We assume from now on
that $R$ is a $\ZZ$-domain of finite type, $\ZZ \subseteq R$, and
that
$f:C \ra \Spec R$ is a smooth projective relative curve. In this
situation the residue class field
$\kappa(\fop)$ of every closed point in $\Spec R$ is a finite field.
We denote by $C_0$ the fiber of the relative curve over the generic point $(0) \in \Spec R$ (in characteristic zero).

For a vector bundle $\shE$ on $C$ we ask when arithmetical Frobenius descent
data imply the semistability of $\shE_0$ on $C_0$. Such data are
given by a sequence $(\fop_n,e_n)_{n \in \NN}$, where $\fop_n$ are
closed points in $\Spec R$ and $e_n$ describes the depth of the
Frobenius descent on $C_{\fop_n}$, i.e., $\shE_{\fop_n} \cong
{F^{e_n}}^* (\shF_n)$, with $\shF_n$ locally free on $C_{\fop_n}$.
We call such a sequence a \emph{Frobenius descent sequence}.

\begin{theorem}
\label{maintheorem}
Let $f:C \ra \Spec R$ be a smooth projective relative curve together with a
fixed $f$-very ample line bundle $\O_C(1)$, where $R$ is a
$\ZZ$-domain of finite type, $\ZZ \subseteq R$. Further let $\shE$ be a locally free
sheaf of rank $r$ on $C$. Suppose there exists a Frobenius descent
sequence $(\fop_n,e_n)_{n \in \NN}$ and assume $(e_n -
|\kappa(\fop_n)|^c)_{n \in \NN} \ra \infty$, where $c$ is the
constant from Lemma \ref{bundlebound}. Then $\shE_0$ is semistable
on the generic fiber $C_0$.
\end{theorem}

\begin{proof}
Let $\shF_n$ be locally free on $C_{\fop_n}$ such that
$\shE_{\fop_n} \cong {F^ {e_n}}^* (\shF_n)$.
First of all we remark that the degree of $\shE$ has to be zero.
This follows, since $\deg(\shE_{\fop_n})= p_n^{e_n} \deg(\shF_n)$
holds for $p_n =
\Char(\kappa(\fop_n))$ and $e_n \ra \infty$ (compare this also with
the proof of \cite[Proposition 3.1]{brennerkaidfrobeniusdescent}).
By Lemma \ref{bundlebound} the number of isomorphism classes of
semistable vector bundles of degree zero and rank $r$ on each fiber
$C_{\fop_n}$ is bounded from above by $|\kappa(\fop_n)|^c$, where
the constant
$c$ depends only on $r$ and the genus of $C$ and is therefore
independent of the closed point $\fop_n$. Since by assumption the
sequence
$(e_n-|\kappa(\fop_n)|^c)_{n \in \NN}$ diverges, we fix
an index $n$ such that
$$e_n-|\kappa(\fop_n)|^c \geq t:= \log_2 (rb)\, ,$$
where $b$ is the bound of Lemma \ref{mumaxbound}. We apply
Lemma
\ref{mumaxboundsemistable} to
$\shE_{\fop_n} \cong {F^{e_n}}^* (\shF_n) = {F^{(e_n-k)}}^*({F^k}^*(\shF_n))$
to conclude that ${F^k}^*( \shF_n)$ is semistable for all $k$ with
$e_n-k > t$.
Due to the inequality $e_n-t \geq |\kappa(\fop_n)|^c$ there has to
be a repetition ${F^k}^* (\shF_n) \cong {F^{k^\prime}}^* (\shF_n)$,
$k < k^\prime < e_n -t$. By Proposition
\ref{equivalentcondiontsstronglysemistable} the locally free sheaf
$\shF_n$ is strongly semistable on $C_{\fop_n}$. Hence
$\shE_{\fop_n}$ is strongly semistable too. By the openness of
semistability we obtain the semistability of $\shE_0 =
\shE|_{C_0}$ on the generic fiber $C_0$.
\end{proof}

\begin{corollary}
In the situation of Theorem \ref{maintheorem} assume that there
exists a closed point $\fop \in \Spec R$ such that $\shE_\fop$
admits an infinite Frobenius descent on $C_\fop$, i.e., we have
$e_n \ra \infty$ in the sequence $(\fop,e_n)_{n \in
\NN}$. Then $\shE_0$ is semistable on the generic fiber $C_0$.
\end{corollary}

\begin{proof}
Since $\fop_n=\fop$ for all $n$, the expression $|\kappa(\fop_n)|^c$ is constant, so this follows from Theorem
\ref{maintheorem}. This corollary can also be proved using
Proposition \ref{equivalentcondiontsflat}, by which the locally free sheaf
$\shE_\fop$ is strongly semistable on $C_\fop$. Again by the
openness of semistability we conclude that $\shE_0$ is semistable on
the generic fiber $C_0$.
\end{proof}

\begin{remark}
In our proof of Theorem \ref{maintheorem} we have found a
closed point $\fop =\fop_n$ in $ \Spec R$ where $\shE_\fop$ is
strongly semistable. Moreover, $\shE_\fop$ is even flat. We have
$\shE_{\fop} \cong {F^{e_n}}^* (\shF_n)$ and
${F^k}^* (\shF_n) \cong {F^{k^\prime}}^* (\shF_n)$,
$k < k^\prime < e_n -t$. Hence
$$\shE_\fop \cong {F^{(e_n-k)}}^*   ({F^{k}}^* (\shF_n) )  \cong {F^{(e_n-k)}}^*   ({F^{k^\prime}}^*
(\shF_n)) \cong {F^{ (k^\prime - k)   }}^* ( {F^{e_n }}^*   (\shF_n))  \, , $$
and $\shE_\fop$ is flat by Proposition
\ref{equivalentcondiontsflat}. A problem of Miyaoka (\cite[Problem 5.4]{miyaokachern})
asks whether for $\shE_0$ semistable the set $S$ of primes
$\fop$ of positive residue characteristic in
$\Spec R$ with strongly semistable reduction $\shE_\fop$ is dense.
It is also open whether for $\shE_0$ semistable there exists always a closed
point $\fop$ such that $\shE_\fop$ is flat (or whether there exists
a Frobenius descent sequence fulfilling the arithmetic condition of
Theorem
\ref{maintheorem}).
\end{remark}

\newcommand{\trans}{{T}}

\begin{example}
Suppose that $D \ra \Spec \FF_p $ is a smooth projective curve and $\shM \neq \O_D$
is an invertible sheaf with $\shM^p \cong \O_D$.
Consider $C = D \times_{\Spec \FF_p} \Spec \FF_p[\trans] \ra \Spec
\FF_p[\trans]$ and $\shL= \shM \times_D C$. For $\fop \in \Spec \FF_p[\trans]$, the sheaf
$\shL_\fop$ is strongly semistable on $C_\fop = D_{\kappa ({\fop})}$. However, it is not
flat at any closed point; this follows from Proposition
\ref{equivalentcondiontsflat}. For every $n$ there exists $q=p^{a_n}$ and an invertible sheaf $\shN$
on
$D_{\FF_q}$ with
$\shL_{\FF_q} \cong \shN^{p^n}$, because this is true on $D_\cloalg {\FF_p}$. Hence we have a Frobenius descent sequence
$(\fop_n,n)$ with $\kappa(\fop_n) = \FF_{p^{a_n}}$.
We have $n - p^{c a_n} \leq 0$ for all $n$, because else there were
a repetition $ \shN^{p^{e'}} =  \shN^{p^{e}}$ for $e < e' \leq n$.
\end{example}

\begin{remark}
The existence of a Frobenius descent for $\shE_\fop$ is by
the so-called Cartier-correspondence \cite[Theorem 5.1]{katznilpotent}
equivalent to the existence of an integrable connection
$\nabla_\fop$ on $\shE_\fop$ with $p$-curvature zero
(see \cite{katznilpotent} and \cite{katzdifferential} for a
detailed exposition of these notions). We neither know whether these
connections are related to each other nor whether they come from a
global connection $\nabla$ on $\shE$. The
Grothendieck-Katz $p$-curvature conjecture
\cite[(I quat)]{katzdifferential} states that given an integrable
connection $\nabla$ on $\shE$ (on $\shX \ra \Spec R$) with the
property that
$\nabla|_{\shE_\fop}$ on $\shX_{\fop}$ has $p$-curvature zero for
$\Char(\kappa(\fop))\gg 0$, then there exists an $\acute{e}$tale
cover $\etcov: Y  \ra \shX_0$ such that
$(\etcov^*(\shE_0),\etcov^*(\nabla_0))$ is trivial. This is only
possible if $\shE_0$ is semistable.
\end{remark}

\section{Counterexamples for the case of constant
characteristic}
\label{counterexamplesection}

In this section we give several examples consisting of rank two vector bundles $\shE$ on
smooth projective relative curves $C$ and Frobenius descent
sequences
$(\fop_n,e_n)_{n \in \NN}$ with certain arithmetical properties, but
where $\shE_0$ is not semistable on the generic fiber
$C_0$.
We first recall briefly our example presented in
\cite{brennerkaidfrobeniusdescent} of a Frobenius descent sequence $(\fop_n, e_n)_{n \in
\NN} =((p_n),1)_{n \in \NN}$ where $p_n$ runs through infinitely many prime numbers. Hence Theorem
\ref{maintheorem} does not hold when the sequence $(e_n)_{n \in
\NN}$ is constant and $\Char \kappa(\fop_n) \to \infty$.
In Example \ref{counterexampleconstantchar} and Example \ref{counterexampleconstantthree}
we present Frobenius descent sequences
of closed points of constant residue characteristic and $e_n \to
\infty$, where again $\shE_0$ is not semistable.
Hence infinite Frobenius descent in different points of the same
residue characteristic does not imply semistability. These examples
are based on examples by P. Monsky showing how the Hilbert-Kunz
multiplicity depends on algebraic parameters, and use the geometric
interpretation of Hilbert-Kunz theory in terms of (strong)
semistability. We do not know of an example where
$|\kappa(\fop_n)|
\to
\infty$ and $e_n \geq |\kappa(\fop_n)|$ and where $\shE_0$ is not
semistable.

\begin{example}
Consider the smooth relative Fermat curve
$$C:= \Proj(\ZZ_d[X,Y,Z]/(X^d+Y^d+Z^d)) \lto
\Spec \ZZ_d$$ of degree $d=2 \ell +1$ and $\ell \geq 2$. The first Frobenius pull-back
of the syzygy bundle $\Syz(X^2,Y^2,Z^2)(3)$ is sitting inside the
short exact sequence
$$0 \lto \O_{C_p}(\ell-1) \lto \Syz(X^{2p},Y^{2p},Z^{2p})(3p) \lto
\O_{C_p}(-\ell+1) \lto 0$$ on the fibers $C_p$ for primes of the
form $p \equiv \ell (d)$ (see \cite[Lemma
2.1]{brennerkaidfrobeniusdescent}). As a consequence of the Cartier-correspondence
this sequence does not split for almost all such $p$ (see \cite[Lemma 2.4] {brennerkaidfrobeniusdescent}).
Therefore, we take the vector bundle $\shE$ on
$C$ defined by the extension corresponding to the
non-trivial $\rm\check{C}$ech-cohomology class $c = Z^{d-1}/XY \in
H^1(C,\O_C(d-3))$. On the fibers $C_p$ we have (because $\dim H^1(C_p,\O_{C_p}(d-3))=1$) the equality
$c= \lambda c^\prime$, where $c^\prime$ is the class corresponding to
our short exact sequence above, hence
$\shE_p \cong F^*(\Syz(X^2,Y^2,Z^2)(3))$. But the
restriction $\shE_0$ to the generic fiber $C_0$ is not semistable.
\end{example}

\begin{lemma}
\label{pullbackhnfiltration}
Let $G:=Z^4+XYZ^2+Z(X^3+Y^3)+ (\algu +\algu^2) X^2Y^2 \in
\FF_2[\algu][X,Y,Z]$ be the equation defining the relative projective curve
$$C:=\Proj(\FF_2[\algu][X,Y,Z]/(G)) \lto \Spec \FF_2[\algu] = \AA^1_{\FF_2}.$$
Further let $\fop \in \Spec \FF_2[\algu]$ be a closed point
corresponding to an algebraic value $\algu \mapsto \alguval \in
\cloalg {\FF_2} $, $\alguval \neq 0,1$, and define
$d:=[\FF_2(\alguval):\FF_2]$. Then the Harder-Narasimhan filtration of
the $d$th Frobenius pull-back of $\Syz(X^2,Y^2,Z^2)(3)$ is given by
the short exact sequence
$$0 \lto \O_{C_\fop}(1) \lto \Syz(X^{2^{d+1}},Y^{2^{d+1}},Z^{2^{d+1}})(3 \cdot
2^d) \lto \O_{C_\fop}(-1) \lto 0$$ on the fiber $C_\fop$.
\end{lemma}

\begin{proof}
To prove the Lemma we use, as in \cite[Lemma
2.1]{brennerkaidfrobeniusdescent}, Hilbert-Kunz theory and its
geometric approach (cf. \cite{brennerhilbertkunz} and
\cite{trivedihilbertkunz}). P. Monsky has shown in \cite[Theorem
4.13]{monskypoints4quartics} that the Hilbert-Kunz multiplicity
(see \cite{monskyhilbertkunz}) of the homogeneous coordinate ring
$R_\fop$ of the fiber $C_{\fop}$ (which is smooth for $\alguval \neq 0,1$) equals
\[e_{HK}(R_\fop) = \begin{cases} 3 \mbox{ if } \kappa(\fop)=\FF_2(\algu),\\
3+\frac{1}{4^d} \mbox{ if } \kappa(\fop)=\FF_2(\alguval) \subset
\cloalg {\FF_2} .
\end{cases} \]
By \cite[Corollary 4.6(ii)]{brennerhilbertkunz} the restriction
$\Omega_{\PP^2}|_{C_\fop} \cong \Syz(X,Y,Z)$ of the cotangent
bundle on the projective plane to the fiber $C_\fop$ is strongly
semistable for the generic point $\fop \in \Spec \FF_2[\algu]$
corresponding to a transcendental value and not strongly semistable
for all algebraic instances. Monsky further proved in \cite[Theorem
3.1]{monskypoints4quartics} that the $d+1$th pull-back
${F^{(d+1)}}^*(\Syz(X,Y,Z))(\frac{3q}{2}-1) \cong
\Syz(X^q,Y^q,Z^q)(\frac{3q}{2}-1)$, $q=2^{d+1}$, has a non-trivial
section. We have to show that this section has no zeros. Consider
the Harder-Narasimhan filtration
$$0 \lto \shL \lto \Syz(X^q,Y^q,Z^q) \lto \tilde{\shL} \lto 0$$
of $\Syz(X^q,Y^q,Z^q)$ where $\shL$ is a line bundle of degree
$\deg(\shL)=-6q+ \alpha$, $\alpha >0$, and the quotient $\tilde{\shL}$ is a line
bundle of degree $\deg(\tilde{\shL})=-6q-\alpha$ (we have to have
$\deg(\shL) + \deg(\tilde{\shL})=-12q=\deg(\Syz(X^q,Y^q,Z^q))$).
By \cite[Corollary 4.6]{brennerhilbertkunz} we can compute
$e_{HK}(R_\fop)$ from the short exact sequence above which
constitutes the strong Harder-Narasimhan filtration in the sense of
\cite[Paragraph 2.6]{langersemistable} and \cite[Section
1]{brennerhilbertkunz} since we are dealing with rank two vector
bundles. We obtain $e_{HK}(R_\fop)=3+\frac{\alpha^2}{4 \cdot
4^{d+1}}=3+\frac{1}{4^d}$ which yields $\alpha = 4$. Hence, $\shL
\otimes_{\O_{C_\fop}} \O_{C_\fop}(\frac{3q}{2}-1)$ has degree zero
and a non-trivial section, so $\shL \cong
\O_{C_\fop}(1-\frac{3q}{2})$ and $\tilde{\shL}\cong
\O_{C_\fop}(-\frac{3q}{2}-1)$. Eventually, we obtain that
$$0 \lto \O_{C_\fop}(1) \lto \Syz(X^q,Y^q,Z^q)(3 \cdot 2^d)\lto \O_{C_\fop}(-1) \lto 0$$
constitutes the Harder-Narasimhan filtration of
$\Syz(X^q,Y^q,Z^q)(\frac{3q}{2})$.
\end{proof}

\newcommand{\denom}{{P(\algu)}}
\begin{example}
\label{counterexampleconstantchar}
In the following we work over the smooth projective relative curve
$$C:=\Proj(\ZZ[\algu]_\denom [X,Y,Z]/(G)) \lto \Spec \ZZ[\algu]_\denom,$$
where $G$ is defined as in Lemma \ref{pullbackhnfiltration} and
where $\denom$ ensures that the relative curve is smooth ($2$ is not
a factor of $\denom$). We consider the vector bundle
$\shE:=\O_C(1)
\oplus  \O_C(-1)$ and the sequence $(\fop_d, e_d)_{d \in
\NN}=((2,f_d(\algu)),d)_{d \in \NN}$ of closed points $\fop_d
\in \Spec \ZZ[\algu]$ of constant residue characteristic two,
where $f_d(\algu)$ denotes an irreducible polynomial of degree
$d$ in $\FF_2[\algu]$. Each fiber $C_{\fop_d}:=C \times_{\Spec
\ZZ[\algu]} \Spec \kappa(\fop_d )$ is a smooth projective curve
defined over $\FF_{2^d}= \kappa(\fop_d)$. The sequence
$$0 \lto \O_{C_{\fop_d}}(1) \lto \Syz(X^{2^{d+1}},Y^{2^{d+1}},Z^{2^{d+1}})(3 \cdot
2^d) \lto \O_{C_{\fop_d}}(-1) \lto 0$$ from Lemma
\ref{pullbackhnfiltration} defines a cohomology class
$$c \in H^1(C_{\fop_d},\O_{C_{\fop_d}}(2)) \cong
\Ext^1(\O_{C_{\fop_d}}(-1),\O_{C_{\fop_d}}(1)).$$ But by
Serre-duality we have $H^1(C_{\fop_d},\O_{C_{\fop_d}}(2)) \cong
H^0(C_{\fop_d}, \O_{C_{\fop_d}}(-1)) = 0$ on each fiber
$C_{\fop_d}$. Hence, the sequence above splits and the class $c$
corresponds to the trivial extension
$\shE_{\fop_d} \cong \O_{C_{\fop_d}}(1)
\oplus
\O_{C_{\fop_d}}(-1)$. Therefore, the vector bundle $\shE$ admits for every closed
point $\fop_d$ a $d$th Frobenius descent $\shE_{\fop_d} \cong
{F^d}^*(\Syz(X^2,Y^2,Z^2)(3))$, but $\shE$ is obviously not
semistable on the generic fiber $C_0$.
Hence $\shE$ has arbitrarily deep Frobenius descent for a sequence of
closed points of constant characteristic, and $e_d \to \infty$. Here
$e_d - |\kappa(\fop_d)|^c = d - 2^{dc} \to - \infty$, and the
numerical condition in Theorem \ref{maintheorem} is not fulfilled.
\end{example}

\begin{lemma}
\label{hnfiltrationpullbackcharthree}
Let $H:=Z^4-XY(X+Y)(X+\algu Y) \in \FF_3[\algu][X,Y,Z]$ and
consider the projective relative curve
$$C:=\Proj(\FF_3[\algu][X,Y,Z]/(H)) \lto \Spec
\FF_3[\algu] = \AA^1_{\FF_3}.$$ Then for every closed point
$\fop \in \Spec \FF_3[\algu]$ corresponding to an algebraic
value $\algu \mapsto \alguval$, $\alguval \neq 0,1$, of degree
$d=[\FF_3(\alguval):\FF_3]$ the
$n=(d+1)$th
Frobenius pull-back of $\Syz(X,Y,Z)$ on $C_\fop$ is sitting inside
the short exact sequence
$$0 \lto \O_{C_\fop}\lto \Syz(X^{3^n},Y^{3^n},Z^{3^n})(\degu) \lto
\O_{C_\fop}(-3) \lto 0$$ where $\degu:=\frac{3^{n+1}-3}{2}$. This sequence splits and
constitutes the Harder-Narasim\-han filtration of
$\Syz(X^{3^n},Y^{3^n},Z^{3^n})(\degu)$.
\end{lemma}

\begin{proof}
It follows from \cite[Proposition 5.10]{trivedihilbertkunz} that
${F^e}^*(\Syz(X,Y,Z))$ is semi\-stable for $e=0 \komdots d-1$ and
not semistable for $e \geq d$. As usual we have
${F^n}^*(\Syz(X,Y,Z)) \cong \Syz(X^{3^n},Y^{3^n},Z^{3^n})$. Let
$$0 \lto \shL \lto \Syz(X^{3^n},Y^{3^n},Z^{3^n}) \lto \tilde{\shL}
\lto 0$$ be the Harder-Narasimhan filtration of the $n$th
pull-back, where $\shL$ and $\tilde{\shL}$ are line bundles of
degrees $\deg(\shL)=k+ \alpha$, where $\alpha >0$, and
$\deg(\tilde{\shL})=k-\alpha$ with $k=-2 \cdot
3^{n+1}=\deg(\Syz(X^{3^n},Y^{3^n},Z^{3^n}))/2$. It was shown by P.
Monsky \cite[Theorem III]{monskyzdp4} that the
Hilbert-Kunz multiplicity of the homogeneous coordinate ring
$R_\fop$ of each fiber $C_\fop$ equals $e_{HK}(R_\fop)=3 +
\frac{1}{9^n}$. Hence we obtain $\alpha = 6$ via \cite[Corollary 4.6]{brennerhilbertkunz}.
Next we tensor the short exact sequence
above with $\degu = \frac{3^{n+1}-3}{2}$ so that $\shL(\degu)$ has
degree zero and $\tilde{\shL}(\degu)$ has degree $-12$. We have to
show that
$\shL(\degu) \cong \O_{C_\fop}$. Assume, $\shL(\degu)$ is not isomorphic
to $\O_{C_\fop}$ (and so $\tilde{\shL}(\degu) \not \cong
\O_{C_\fop} (-3))$. We recall that one can compute from the presenting
sequence of
$\Syz(X^{3^n},Y^{3^n},Z^{3^n})(m)$ the
Hilbert-Kunz function via
\begin{eqnarray*}
\length((R_\fop/(X,Y,Z)^{[3^n]})_m)&=&h^0(C_\fop ,\O_{C_\fop}(m))-3
h^0(C_\fop
,\O_{C_\fop}(m-3^n))\\
& & + h^0(C_\fop ,\Syz(X^{3^n},Y^{3^n},Z^{3^n})(m)) \, ,
\end{eqnarray*}
and summing these expressions over $m \geq 0$ (where we can restrict
to a finite sum, as the alternating sum is $0$ for $m \gg 0$).
In fact we only have to sum up to $\degu+3$, since by
Serre-duality we have
$H^1(C_\fop ,\shL(\degu+4)) \cong H^0(C_\fop ,(\shL(\degu))^*(-3))^*=0$ and
$H^1(C_\fop ,\tilde{\shL}(\degu+4)) \cong H^0(C_\fop ,(\tilde{\shL}(\degu))^*(-3))^*=0$
and therefore the cohomology
$H^1(C_\fop ,\Syz(X^{3^n},Y^{3^n},Z^{3^n})(\degu+4))$ vanishes. For $m=\degu+1,\degu+2 ,\degu+3$
we get by Riemann-Roch (the quotient line bundle
$\tilde{\shL}(m)$ has no non-trivial global sections in these
degrees):
\begin{eqnarray*}
h^0(C_\fop ,\Syz(X^{3^n},Y^{3^n},Z^{3^n})(\degu+1))\!\!&=&\!\!h^0(C_\fop ,\shL(\degu+1))\,\,=\,\,2,\\
h^0(C_\fop ,\Syz(X^{3^n},Y^{3^n},Z^{3^n})(\degu+2))\!\! &=&\!\!h^0(C_\fop ,\shL(\degu+2))\,\,=\,\,6,\\
h^0(C_\fop ,\Syz(X^{3^n},Y^{3^n},Z^{3^n})(\degu+3))\!\!&=&\!\!
h^0(C_\fop,\shL(\degu+3))\,\,=\,\,10.
\end{eqnarray*}
In degrees $m \leq \degu$ we have
$h^0(C_\fop ,\Syz(X^{3^n},Y^{3^n},Z^{3^n})(m))=0$ (for $m = \degu$ by assumption). The Hilbert function
$HF_{R_\fop} (m)= h^0(C_\fop, \O_{C_\fop}(m))$
of the homogeneous coordinate ring $R_\fop$ of each fiber $C_\fop$
is given by $HF_{R_\fop}(0)=1$, $HF_{R_\fop}(1)=3$ and
$HF_{R_\fop}(m)=4m-2$ for $m \geq 2$. So we have accumulated all data
to compute the Hilbert-Kunz function $\varphi$ for the exponent
$n=d+1$ and get:
\begin{eqnarray*}
\varphi(3^n) &=& \sum_{m=0}^{\degu+3} h^0(C_\fop ,\O_{C_\fop}(m)) - 3
\sum_{m=0}^{\degu+3} h^0(C_\fop ,\O_{C_\fop}(m-3^n)) +18\\
&=& 4 + \sum_{m=2}^{\degu+3}(4m-2) -3(4 + \sum_{m=2}^{\degu+3-3^n}(4m-2))
+ 18\\
&=& 2(\degu+3)^2-2-3(2(\degu+3-3^n)^2-2)+10\\
&=&-4(\degu+3)^2+4(\degu+3)3^{n+1}-2 \cdot 3^{2n+1}+14\\
&=&-(3^{n+1}+3)^2+2(3^{n+1}+3)3^{n+1}-2 \cdot 3^{2n+1} +14\\
&=&3^{2n+1} + 5.
\end{eqnarray*}
But this contradicts \cite[Theorem IIb]{monskyzdp4}, where P.
Monsky proved with different methods that $\varphi(3^n)=3^{2n+1}+9$. Hence,
$\shL(\degu) \cong \O_{C_\fop}$ and in particular $\tilde{\shL}(\degu)
\cong \O_{C_\fop}(-3)$.
The sequence splits since $H^1(C_\fop, \O_{C_\fop}(3))=0$.
\end{proof}

\begin{example}
\label{counterexampleconstantthree}
First we look at the smooth projective relative curve
$$C:=\Proj(\ZZ[\algu]_\denom [X,Y,Z]/(H)) \lto \Spec \ZZ[\algu]_\denom,$$
where $H$ is defined as in Lemma
\ref{hnfiltrationpullbackcharthree} and localization at $\denom$ makes the curve smooth
($3$ is not a factor of $\denom$).
We consider a sequence
$(\fop_d)_{d \in \NN} = ((3,f_d(\algu)))_{d \in \NN}$ of closed
points in $\Spec \ZZ[\algu]$, where $f_d(\algu)$ denotes a
polynomial of degree $d$ which is irreducible in
$\FF_3[\algu]$.
Unlike in Example \ref{counterexampleconstantchar} no twist of a
pull-back of $\Syz(X,Y,Z)$ on a fiber $C_{\fop_d}$ has degree zero.
Therefore, we pass over to a suitable finite cover of $C$.
As in \cite[Remark 4.4]{brennermonskytightclosure} we look at the
ring homomorphism
$$\ZZ[\algu][X,Y,Z]/(H) \lto
\ZZ[\algu][X,Y,Z]/(W^8 - (U^4+V^4)(U^4+\algu V^4)) =: B$$ given
by $X \mapsto U^4$, $Y \mapsto V^4$, $Z \mapsto UVW^2$ and the
corresponding map $g:D:=\Proj(B) \rightarrow C$ of smooth projective
relative curves (over $\Spec \ZZ[\algu]_\denom$). Let $\shE:=
\O_D(6)
\oplus \O_D(-6)$. It is easy to see that
$g^*(\O_C(1)) \cong \O_D(4)$ and hence the pull-back under $g$ of
the short exact sequence of Lemma
\ref{hnfiltrationpullbackcharthree} yields
$$0 \lto \O_{D_{\fop_d}} \lto \Syz(U^{4\cdot 3^n},V^{4 \cdot 3^n},U^{3^n}V^{3^n}W^{2\cdot 3^n})(4u) \lto \O_{D_{\fop_d}}(-12) \lto
0$$
($n=d+1$) on each fiber $D_{\fop_d}$. After tensoring with $\O_{D_{\fop_d}}(6)$ we obtain
$$0 \lto \O_{D_{\fop_d}}(6) \lto {F^{(d+1)}}^*(\Syz(U^4,V^4,UVW^2)(6)) \lto \O_{D_{\fop_d}}(-6) \lto
0\, .$$
Because this sequence splits we get
$\shE_{\fop_d} \cong {F^{(d+1)}}^*(\Syz(U^4,V^4,UVW^2)(6))$. Therefore $(\fop_d, d+1)_{d \in \NN}$ is a Frobenius descent sequence.
\end{example}

\section*{Acknowledgements}

We thank the referee and H\'{e}l\`{e}ne Esnault for many useful comments and Helena Fischbacher-Weitz for careful reading of this manuscript.

\bibliographystyle{mrl}

\end{document}